\documentclass[12pt,a4paper]{article}
\usepackage{amsmath,amsthm,amssymb,amscd}
\usepackage{latexsym}
\usepackage{indentfirst}
\usepackage{bibentry}
\usepackage{textcomp}
\usepackage{float}
\usepackage{multirow}
\usepackage{booktabs}
\usepackage{graphicx}
\usepackage{color}
\usepackage{cite}
\usepackage{authblk}
\usepackage{algorithmicx,algorithm}
\usepackage{listings}
\usepackage{xcolor}
\usepackage{alltt}
\usepackage{enumitem}
\usepackage{lscape}
\usepackage{geometry}
\usepackage{enumitem}
\usepackage{longtable}
\newcommand*{\thead}[1]{\multicolumn{1}{c}{\bfseries #1}}

\usepackage{afterpage}
\usepackage{capt-of}
\makeatletter \@addtoreset{equation}{section}

\makeatother

\theoremstyle{definition}

\begin{document}

\title{Smallest Eigenvalue of Large Hankel Matrices at Critical Point: Comparing Conjecture With Parallelised Computation}

\author[a]{{Yang Chen}\footnote{yangbrookchen@yahoo.co.uk,~yayangchen@umac.mo}}
\author[a]{{Jakub Sikorowski}\footnote{sikorowski@gmail.com}}
\author[b]{{Mengkun Zhu}\footnote{Corresponding author: Zhu\_mengkun@163.com}}

\affil[a]{Department of Mathematics, University of Macau,
Avenida da Universidade, Taipa, Macau, China}
\affil[b]{School of Science, Qilu University of Technology (Shandong Academy of Sciences)\\
Jinan 250353, China}

\renewcommand\Authands{ and }

\maketitle

\begin{abstract}

We propose a novel parallel numerical algorithm for calculating the smallest eigenvalues of highly ill-conditioned matrices. It is based on the {\it LDLT} decomposition and involves finding a $k \times k$ sub-matrix of the inverse of the original $N \times N$ Hankel matrix $H_N^{-1}$ . The computation involves
extremely high precision arithmetic, message passing interface, and shared memory parallelisation. We demonstrate that this approach achieves good
scalability on a high performance computing cluster (HPCC) which constitute a major improvement of the earlier approaches.
We use this method to study a family of Hankel matrices generated by the weight $w(x)={\rm e}^{-x^{\beta}},$ supported on $[0,\infty)$ and $\beta>0.$
Such weight generates Hankel determinant, a fundamental object in random matrix theory. In the situation where $\beta>1/2,$  the smallest eigenvalue
tend to 0, exponentially fast as $N$ gets large. If $\beta<1/2,$ the situation where the classical moment problem is indeterminate,
the smallest eigenvalue is bounded from below by a positive number for all $N$, including infinity. If $\beta=1/2,$ it is conjectured that the smallest eigenvalue tends to 0
algebraically, with a precise exponent. The algorithm run on the HPCC producing fantastic match between the theoretical value of $2/\pi$ and the numerical result.
\end{abstract}

Keywords:  Random matrix, Parallel eigensolver, Extremely ill-conditioned Hankel matrices, Smallest eigenvalue
\\

MSC: 15B57, 34E05, 42C05, 65F15, 65Y05

\section{Background and motivation}
\noindent Hankel matrices has entries that are moments \cite{C16,new11} of probabilities measures or weight functions. And plays an important role in the theory of Hermitian random matrices \cite{C17}. The study of the largest and smallest eigenvalues are important since they provide useful information about the nature of the Hankel matrix generated by a given density, e.g. they are related with the inversion of Hankel matrices, where the condition numbers are enormously large.

Given $\left\{\mu_{j}\right\}$ the moment sequence of a weight function $w(x)(>0)$ with infinite support $s$,
\begin{equation}\label{b2}
\mu_{j}:=\int_{s} x^{j}w(x)dx,~~k=0,1,2,\ldots,
\end{equation}
the Hankel matrices, it is known that
\begin{equation*}
H_{N}:=\left(\mu_{i+j}\right)_{i,j=0}^{N-1},~~ N=1,2,\ldots
\end{equation*}
are positive definite \cite{zzz}.

Let $\lambda_{1}$ be the smallest eigenvalue of $H_{N}$. The asymptotic behavior of $\lambda_{1}$ for large $N$ has broad interest, see e.g. \cite{CL1,BCI,ECW,C14,C10,C3,C5,C6,new33,C22}. The authors in \cite{C24,C25} have studied the behavior of the condition number ${\rm cond}\left(H_{N}\right):=\frac{\lambda_{N}}{\lambda_{1}}$, where $\lambda_{N}$ denotes the largest eigenvalue of $H_{N}$.

Szeg\"{o}{\color{black}\cite{C3}} investigated the asymptotic behavior of $\lambda_{1}$ for the Hermite weight $w(x)={\rm e}^{-x^{2}},x\in\mathbb{R}$ and the classical Laguerre weight $w(x)={\rm e}^{-x},x\geq0$. He found\footnote[1]{In all of this paper, $a_{N}\simeq b_{N}$ means $\lim_{N\rightarrow\infty}a_{N}/b_{N}$=1.}
\begin{equation*}
\lambda_{1}\simeq A N^{\frac{1}{4}}B^{\sqrt{N}},
\end{equation*}
where $A,B$ are certain constants, satisfying $0<A,~0<B<1$. Moreover, in the same paper, it was showed that the largest eigenvalue $\lambda_{N}$ refers to the Hankel matrices $\left[\frac{1}{i+j+1}\right]_{i,j=0}^{N-1}$, $\left[\Gamma\left(\frac{i+j+1}{2}\right)\right]_{i,j=0}^{N-1}$ and $\left[\Gamma(i+j+1)\right]_{i,j=0}^{N-1}$ were approximated by $\frac{\pi}{2}$, $\Gamma\left(N+\frac{1}{2}\right)$ and $(2N)!$ respectively. Widom and Wilf \cite{C6} studied the situation that the density $w(x)$ is supported on a compact interval $[a,b]$, such that the Szeg\"{o} condition
\begin{equation*}
\int_{a}^{b}\frac{-\ln w(x)}{\sqrt{(b-x)(x-a)}}dx<\infty,
\end{equation*}
holds, they found that
\begin{equation*}
\lambda_{1}\simeq A \sqrt{N}B^{N}.
\end{equation*}

Chen and Lawrence \cite{CL1} found the asymptotic behavior of $\lambda_{1}$ with the weight function $w(x)={\rm e}^{-x^{\beta}},~x\in[0,\infty),~\beta>\frac{1}{2}$. Critical point at $\beta=\frac{1}{2}$. It marks the transition point at which the moment problem becomes indeterminate. Theoretical challenges to finding the asymptotic behavior at critical point. Berg, Chen and Ismail \cite{BCI} proved that the moment sequence (\ref{b2}) is determinate iff $\lambda_{1}\rightarrow0$ as $N\rightarrow\infty$. This is a new criteria for the determinacy of the Hamburger moment problem. Chen and Lubinsky \cite{C10} found the behavior of $\lambda_{1}$ when $w(x)={\rm e}^{-|x|^{\alpha}},~x\in\mathbb{R},~\alpha>1$. Berg and Szwarc \cite{C14} proved that $\lambda_{1}$ has exponential decay to zero for any measure which has compact support. Recently, Zhu \emph{et al.} \cite{C22} studied the Jacobi weight, $w(x)=x^{\alpha}(1-x)^{\beta},~x\in[0,1],~\alpha>-1,~\beta>-1$ and derived a approximation formula of $\lambda_{1}$,
\begin{equation*}
\lambda_{1}\simeq2^{\frac{15}{4}}\pi^{\frac{3}{2}}\big(1+2^{\frac{1}{2}}\big)^{-2\alpha}\big(1+2^{-\frac{1}{2}}\big)^{-2\beta}N^{\frac{1}{2}}\big(1+2^{\frac{1}{2}}\big)^{-4(N+1)},
\end{equation*}
which reduces to Sezg\"{o}'s result \cite{C3} if $\alpha=\beta=0$.

This paper is concerned with a numerical computation that is motivated by random matrix theory. We are interested in finding the numerical value of the smallest eigenvalue of a class of $\mathbb{R}^{N\times N}$ Hankel matrices.

\section{Mathematical statement of the problem}\label{sec:problem}
\noindent We consider the weight function
\begin{equation}\label{weight}
 w(x):= {\rm e}^{-x^\beta}, ~~~~x\in[0,\infty),~~ \beta>0,
 \end{equation}
and the moments, based on (\ref{b2}), are given by
\begin{equation}\label{moments}
\mu_k= \dfrac{1}{\beta}\, \Gamma\left(\dfrac{k+1}{\beta}\right),~~  k = 0,1,2...
\end{equation}
In this case, the entries of Hankel matrix
\begin{align}\label{eq:matrix}
H_N =
\begin{bmatrix}
 \mu_{0}   & \mu_{1} & \mu_{2} &  \ldots    & \ldots     & \mu_{N-1} \\
 \mu_{1}   & \mu_2   &         &            &            & \vdots    \\
 \mu_{2}   &         &  \ddots &            &          & \vdots    \\
 \vdots    &         &         &   \ddots   &            & \mu_{2N-4}\\
 \vdots    &         &         &            & \mu_{2N-4} & \mu_{2N-3}\\
 \mu_{N-1} & \ldots  & \ldots  & \mu_{2N-4} & \mu_{2N-3} & \mu_{2N-2}
\end{bmatrix}
\end{align}
increase factorially along the diagonal. In this paper, we study the smallest eigenvalue at the critical point $\beta=\frac{1}{2}$.

\subsection{Condition number}
Let $\left\{\lambda_{i}\right\}_{i=1}^{N}$ be the $N$ positive eigenvalues of $H_{N}$. In order to highlight the numerical challenge, we would like to estimate the condition number of the problem
\begin{align*}
   {\rm cond}(H_N) &= \lambda_N / \lambda_1,
\end{align*}
where $\lambda_N$ is the largest eigenvalue and $\lambda_1$ is the smallest. We observe that the largest eigenvalues $\lambda_N$ is bounded by the trace of the matrix
\begin{align*}
 \lambda_N \leq \sum_{i=1}^N \lambda_i = \text{tr}(H_N)= \sum_{k=0}^{N-1} \frac{1}{\beta}\Gamma \left(\frac{2k+1}{\beta} \right).
\end{align*}
On the other hand it must be larger or equal than the largest entry along the diagonal, i.e.
\begin{align*}
 \lambda_N \ge \frac{1}{\beta}\Gamma \left(\frac{2N-1}{\beta} \right).
\end{align*}
For sufficiently large $N$, we note
\begin{align*}
 \sum_{k=0}^{N-1} \frac{1}{\beta}\Gamma \left(\frac{2k+1}{\beta}\right)\simeq \frac{1}{\beta}\Gamma \left(\frac{2N-1}{\beta}\right),
\end{align*}
therefore,
\begin{align*}
 \lambda_N \simeq \frac{1}{\beta}\Gamma \left(\frac{2N-1}{\beta} \right).
\end{align*}
For the case of $\beta=\frac{1}{2}$, the smallest eigenvalues are presented in the Table \ref{fig:table2}. They are the order of $\frac{1}{10}$ and slowly decreasing with $N$. Therefore, we estimate the condition number (assuming $\lambda_1=\frac{1}{10}$)
\begin{align*}
{\rm cond}(H_N) &\simeq 20\;\, \Gamma(4N-2).
\end{align*}
 Notice that it grows factorially with increasing $N$, as a result, we say that these matrices are extremely ill-conditioned. This property makes the numerical computation particularly challenging.

The condition number of a problem quantifies the sensitivity of the output value to small changes in the input. For example, in our case, a finite precision of numerical representation of the Hankel matrix introduces a (hopefully small) change of the computed eigenvalue. We want this change to be less than $10^{-15}$. However, an astronomically large condition number makes it hard, and the initial Hankel matrix must be initialized with extreme amount of precision for example $N=4500$ needs $\sim70,000$ digits of precision. Subsequently the intermediate arithmetic operations must be performed with even higher precision in order for the rounding error not to introduce large errors in the output. This extreme precision requirements demand significant computing resources, and presents novel challenges and trade offs.

We would like to stress that the algorithm we use in this paper was chosen to be as numerically stable as possible (see comparison between Secant, Householder, Jacobi, Lanczos in \cite{ECW} and the Section \ref{sec:compare} on page \pageref{sec:compare}), and the numerical challenges are intrinsic to the problem and are not due to instabilities of the algorithms.

\section{Properties of Hankel matrices}\label{sec:properties}
\noindent The Hankel matrices we defined in section \ref{sec:problem} have many interesting properties, here we would like to describe those properties that are used in this paper. We notice that the matrix in eq. (\ref{eq:matrix}) is explicitly symmetric. If $\pi_{N}(x)$ is a polynomial of degree $N$ with real coefficients $c_{j}$, $j=0,1,...N$, namely,
\begin{equation*}
\pi_{N}(x):=\sum_{j=0}^{N}c_{j}x^{j},
\end{equation*}
then the quadratic form
\begin{align*}
 \sum\limits_{i,\,j=0}^{N} c_i \mu_{i+j}c_j= \int_0^{\infty}\left[\pi_{N}(x)\right]^2 w(x) dx,
\end{align*}
is positive definite.

The authors in \cite{CL1} have studied the small eigenvalues with respect to the weight function (\ref{weight}) and derived the asymptotic expression for smallest eigenvalue for $\beta>\frac{1}{2}$. They have found that the smallest eigenvalue decreases exponentially with increasing $N$. On the other hand, the situation for $0<\beta<\frac{1}{2}$, the weight (\ref{weight}) is Stieltjes indeterminate\footnote{There are infinitely many measures support on $[0,\infty)$ with the moments (\ref{moments}).}, they have found that the lower bound of the smallest eigenvalue is greater than a positive number for any $N$.

At the critical point $\beta=\frac{1}{2}$, the asymptotic behavior changes abruptly, and we have a phase transition. Furthermore, the authors of \cite{CL1} argued that under certain assumptions the asymptotic expression for the smallest eigenvalue is
\begin{align}\label{eq:leading_order_old}
 \lambda_1(N) &\simeq 8\pi \dfrac{\sqrt{\log(4\pi N {\rm e})}}{(4\pi N {\rm e})^{\frac{2}{\pi}}},~~~ \text{ for } \beta= \dfrac{1}{2},
\end{align}
for large $N$. Some of those assumptions could not be proven, therefore, we treat eq. (\ref{eq:leading_order_old}) as a conjecture. In the following subsections \ref{sec:pol} and \ref{sec:sad} we summarise main steps of \cite{CL1} performed to find eq. (\ref{eq:leading_order_old}).

\subsection{Polynomials $P_n$}\label{sec:pol}
The main heroes in the derivation of eq. (\ref{eq:leading_order_old}) are $P_j(x)$, the orthonormal polynomials associated with $w(x)$
\begin{align*}
 \int_0^{\infty} P_i(x)  P_j(x){\rm e}^{-x^{\beta}} dx = \delta_{ij},
\end{align*}
where $\delta_{ij}$ denotes the Kronecker's delta, i.e., $\delta_{ij}=1$ if $i=j$ and $\delta_{ij}=0$ for $i\neq j$.

It is shown \cite{CL1} that the smallest eigenvalue $\lambda_1(N)$ is bounded from below
\begin{align}\label{eq:eig_approx}
 \lambda_1(N)\ge \dfrac{2\pi}{\sum\nolimits_{j=0}^{N} K_{jj}},
\end{align}
where
\begin{equation*}
K_{jk}:= \int_{-\pi}^{\pi}P_j\left({\rm -e}^{{\rm i}\phi}\right)P_k\left({\rm -e}^{-{\rm i}\phi}\right)d\phi.
\end{equation*}

Using the Christoffel-Darboux formula \cite{zzz} and the result presented in \cite{C23} for large $N$ off-diagonal recurrence coefficients, it was found that
\begin{align}\label{eq:K_in_P}
 \sum_{j=0}^N K_{jj} \simeq \pi^2 N^2 \int_{-\pi}^{\pi} \frac{P_N\left(-{\rm e}^{{\rm i}\phi}\right)P_{N+1}\left(-{\rm e}^{-{\rm i}\phi}\right) - P_N\left(-{\rm e}^{-{\rm i}\phi}\right) P_{N+1}\left(-{\rm e}^{{\rm i}\phi}\right)}{{\rm e}^{{\rm i}\phi}-{\rm e}^{-{\rm i}\phi}} d\phi
\end{align}
Finally, the authors of \cite{CL1} used the asymptotic expression for $P_N(x)$ for large $N$
\begin{align}\label{eq:poly}
 P_N(x) \simeq\dfrac{(-1)^N}{2\pi}(-x)^{-\frac{1}{4}} N^{-\frac{1}{2}}
         \exp\left[
            \dfrac{\sqrt{-x}}{\pi}
             \log
                \left(
                  \dfrac{4\pi N {\rm e}}{\sqrt{-x}}
                \right)
            \right],
\end{align}
to evaluate eq. (\ref{eq:K_in_P}) for large $N$. We perform the integration in the following section.

\subsection{Saddle point approximation}\label{sec:sad}
In this section we evaluate the integral in eq. (\ref{eq:K_in_P}) using the Laplace method, see \cite{C27}. We start by substituting the asymptotic expressions eq. (\ref{eq:poly}) into eq. (\ref{eq:K_in_P}). This leads to an integral with an integrand that has a maximum at $\phi=0$ and falls off to $0$ towards $\phi=\pm\pi$
\begin{align}
 \sum _{j=0}^N K_{jj} \simeq - \int_{-\pi}^{\pi} &d\phi \dfrac{N}{4 \left({\rm e}^{{\rm i}\phi}-{\rm e}^{-{\rm i}\phi}\right)}\nonumber\\
   \Bigg\{
         &\exp\left[
            \frac{{\rm e}^{ \frac{{\rm i}\phi}{2}}}{\pi}
              \log
                \Big(
                  4\pi N    {\rm e}^{1-\frac{{\rm i}\phi}{2}}
                \Big)
           +\frac{{\rm e}^{-\frac{{\rm i}\phi}{2}}}{\pi}
              \log
                \Big(
                  4\pi(N+1) {\rm e}^{1+\frac{{\rm i}\phi}{2}}
                \Big)
                  \right]\nonumber\\
        -&\exp\left[
            \frac{{\rm e}^{-\frac{{\rm i}\phi}{2}}}{\pi}
              \log
                \Big(
                  4\pi N {\rm e}^{1+\frac{{\rm i}\phi}{2}}
                \Big)
           +\frac{{\rm e}^{ \frac{{\rm i}\phi}{2}}}{\pi}
                \log
                \Big(
                  4\pi(N+1) {\rm e}^{1-\frac{{\rm i}\phi}{2}}
                \Big)
                  \right]\Bigg\}
\end{align}
Furthermore, the width of the peak decreases as $N$ increases, as a result, we can use the Laplace method to evaluate the integral for large $N$. We recast the integrand as ${\rm e}^{-f(\phi)}$ and expand $f(\phi)$ in Taylor series
\begin{align}
 \sum\limits_{j=0}^N K_{jj} \simeq
 \int_{-\pi}^{\pi} {\rm e}^{-f(\phi)}d\phi \simeq
 \int_{-\infty}^{\infty} {\rm e}^{-                 f_0
                                   -                 f_1\; \phi
                                   - \frac{1}{2 }\; f_2\; \phi^2
                                   - \frac{1}{6 }\; f_3\; \phi^3
                                   - \frac{1}{24}\; f_4\; \phi^4 +...} d\phi
\end{align}
 Expanding around $\phi=0$ we find that $f_{1}$ and $f_3$ vanish. Moreover,
 \begin{align*}
  \frac{f_4}{f_2^2} = \mathcal{O}\left(\frac{1}{\log(N)}\right) \ll 1
 \end{align*}
for large $N$, therefore, we can treat $f_4$ and higher order terms as small perturbation\footnote{Similarly $\frac{f_{2n}}{f_2^n} =\mathcal{O}\left(\frac{1}{\log^n(N)}\right)$ and $f_{2n+1}=0$ for $n=1,2,3...$} and expand the exponential
\begin{align}\label{eq:2nd_order}
 \sum\limits_{j=0}^N K_{jj} &\simeq
   {\rm e}^{-f_0} \sqrt{\dfrac{2\pi}{f_2}}\: \left(1 - \dfrac{f_4}{8\,f_2^2}+... \right),\\
\end{align}

 where
 \begin{align*}
   {\rm e}^{-f_0} &= 2^{-3 +\frac{4}{\pi}}\, {\rm e}^{\frac{2}{\pi}}\,
             \pi^{\frac{2}{\pi}-1}\,  N^{\frac{3}{2}+\frac{1}{\pi}}\,
             (N+1)^{-\frac{1}{2}+\frac{1}{\pi}}
             \log\left(1+1/N\right),\nonumber\\
   12\pi^2f_2 &= -3\pi(2+\pi-2\log4\pi)
                    +\log^2 N
                    +3\pi \log(N+1)\nonumber\\
                 &\quad  +\left[3\pi-2\log(N+1)\right]\log N
                    +\log^2(N+1)\nonumber\\
   240\,\pi^4\,f_4 &=
          -30\pi^3 \left(-3+\pi+\log4\pi\right)
          +2\log^4N-15\pi^3\log(N+1)\nonumber\\
   &\quad -8\log^3N\log(N+1)
          -20\pi^2\log^2(N+1)
          +2\log^2(N+1)\nonumber\\
   &\quad -4\log^2N\left[5\pi^2-3\log^2(N+1)\right]\nonumber\\
   &\quad +\left[-15\pi^3+40\pi^2\log(N+1)-8\log^3(N+1) \right]\log N,\nonumber
\end{align*}
Expanding $f_1,...,f_4$ around $N\to\infty$ we obtain
\begin{align*}
   {\rm e}^{f_0} &= 2^{\frac{4}{\pi}-3}
              {\rm e}^{\frac{2}{\pi}}
            \pi^{\frac{2}{\pi}-1}
              N^{\frac{2}{\pi}} \left[1 + \frac{\pi-1}{\pi} \frac{1}{N} + \mathcal{O}\left(\frac{1}{N^2}\right)\right],\nonumber\\
   f_2 &= \frac{1}{2\pi} \left[\log\left(4\pi {\rm e}^{-1-\frac{\pi}{2}} N\right)
        + \frac{1}{4\pi N}
        + \mathcal{O}\left(\frac{1}{N^2}\right)\right]\nonumber\\
   f_4 &=-\frac{1}{8\pi}\left[ \log(4\pi {\rm e}^{\pi-3} N) +\mathcal{O}\left(\frac{1}{N}\right)\right]\nonumber
\end{align*}
Taking only the first term in the above equations and neglecting all corrections to the saddle point approximation in eq. (\ref{eq:2nd_order})
we obtain
\begin{align*}
 \sum\limits_{j=0}^N K_{jj} \simeq \dfrac{(4\pi N {\rm e})^{\frac{2}{\pi}}}
                                          {4\sqrt{\log\left(4\pi {\rm e}^{-1-\frac{\pi}{2}}N\right)}},
\end{align*}
which together with eq. (\ref{eq:eig_approx}) gives\footnote{Authors of this publication have noticed that eq. (\ref{eq:saddle_point}) does not exactly match the eq. (\ref{eq:leading_order_old}) which is quoted from \cite{CL1}, we suspect there was a typo in the earlier publication.}
\begin{align}\label{eq:saddle_point}
 \lambda_1(N) \simeq \dfrac{8\pi \sqrt{\log\left(4\pi {\rm e}^{-1-\frac{\pi}{2}} N\right)}}{(4\pi {\rm e})^{\frac{2}{\pi}}N^{\frac{2}{\pi}}} .
\end{align}

Subsequently including the next-to-leading order correction in the saddle point, i.e. including the $\frac{f_4}{8\,f_2^2}$ term, see eq. (\ref{eq:2nd_order}), we obtain
\begin{align}\label{eq:nlo_order}
 \sum\limits_{j=0}^N K_{jj} &\simeq \dfrac{(4\pi N {\rm e})^{\frac{2}{\pi}}}
                                          {4\sqrt{\log\left(4\pi {\rm e}^{-1-\frac{\pi}{2}}N\right)}}
                                          \left(1 + \dfrac{\pi}{16\log N}\right),\\
 \lambda_1(N) &\simeq \dfrac{8\pi \sqrt{\log\left(4\pi {\rm e}^{-1-\frac{\pi}{2}} N\right)}}{(4\pi {\rm e})^{\frac{2}{\pi}}N^{\frac{2}{\pi}}}
                       \left(1 - \dfrac{\pi}{16\log N}\right)\\
              &\simeq \dfrac{8\pi \sqrt{\log\left(4\pi {\rm e}^{-1-\frac{\pi}{2}-\frac{\pi}{8}} N\right)}}{(4\pi {\rm e})^{\frac{2}{\pi}}N^{\frac{2}{\pi}}} \label{eq:nlo_order_all_e}\\
              &\simeq \dfrac{8\pi\sqrt{\log{N}}}
                                          {(4\pi N {\rm e})^{\frac{2}{\pi}}}
                                          \left[1 + \dfrac{8 \log\left(4\,\pi\, {\rm e}^{-1-\frac{\pi}{2}}\right) -\pi}{16\log N}\right].\label{eq:nlo_order_end}
\end{align}
The transformations from eq. (\ref{eq:nlo_order}) to eq. (\ref{eq:nlo_order_end}) neglect all next-to-next-to-leading order terms. Moreover, we need to stress that this expression might not include all of the next-to-leading order contributions. For example, it does not include next-to-leading order contributions (potentially) coming from
\begin{itemize}
 \item The approximation in eq. (\ref{eq:K_in_P}).
 \item The asymptotic nature of eq.(\ref{eq:poly}).
\end{itemize}
However, we notice that the subleading term we calculated corrects only the log term in eq. (\ref{eq:nlo_order_all_e}) and does not affect the overall $N^{-\frac{2}{\pi}}$ factor, and $-\frac{2}{\pi}$ which we refer to as the leading exponent. In section \ref{sec:exponent} we show that directly computing the determinant of Hankel matrices we were able to find $-\frac{2}{\pi}$ to a good precision.

Keeping in mind that eq. (\ref{eq:nlo_order}) might not capture the full next-to-leading order contribution it would be prudent to use only the leading order approximation
\begin{align}\label{eq:leading_order}
 \sum\limits_{j=0}^N K_{jj} &\simeq \dfrac{(4\pi N {\rm e})^{\frac{2}{\pi}}}
                                           {4\sqrt{\log{N}}},\\
 \lambda_1(N)               &\simeq \dfrac{8\pi\sqrt{\log{N}}}
                                           {(4\pi N {\rm e})^{\frac{2}{\pi}}}.
\end{align}

Notice that the asymptotic expression in eq. (\ref{eq:leading_order_old}) does not decrease exponentially with $N$ as opposed to the $\beta>\frac{1}{2}$ case\cite{CL1}. Moreover, the subleading terms in eq. (\ref{eq:nlo_order_end}) are suppressed only by $\log N$ terms and decrease very slowly with increasing $N$. This makes the difference between eq. (\ref{eq:leading_order}) and the numerically computed also decrease very slowly with $N$. This was indeed observed in \cite{ECW}. The differences were decreasing with $N$, however, much slower than for $\beta\neq1/2$, and, as a result, the numerics did not convincingly confirm eq. (\ref{eq:leading_order_old}). In this paper, we endeavour to improve on that.

\section{Numerical algorithm}\label{sec:methods}
In this section, we describe the novel numerical algorithm for computing the smallest eigenvalues of a highly-ill conditioned matrix. We tweak and optimise the algorithm for Hankel matrix in eq. (\ref{eq:matrix}).

Large parts of the new algorithm are based on the Secant algorithm described in \cite{ECW}. However, for completeness, we will explain the new algorithm without assuming any reader's knowledge of \cite{ECW}.

\subsection{Precision}\label{sec:precision}
The posed problem is a highly ill-condition, therefore, in order to obtain accurate results we have to perform arithmetic operations and store numbers with extremely high precision. In order to obtain appropriate precision, we have used the arbitrary precision integer arithmetic implemented in GNU Multiple Precision Arithmetic Library. Each element of the initial matrix ${\tt H}$ is represented by
\begin{align}\label{eq:precision}
 \dfrac{\mathbb{Z}_{GMP}}{2^K}
\end{align}
where $K$ represents the number of bits of precision of the calculation, and $\mathbb{Z}_{GMP}$ is the GMP arbitrary precision integer. This way we obtain a representation of a number with a fixed number of bits to the right of the decimal point and an arbitrary number of bits to the left of the decimal point. Furthermore, all intermediate computations during the {\it LDLT} decomposition are done with twice as much precision to the right of the decimal point. After, the {\it LDLT} decomposition we truncate the numbers to $K$ bits of precision.

In order to test if the number of bits of precision $K$ was enough to obtain the correct output, we perform the same computation with an increased number of bits of precision. Only when the eigenvalues calculated match the eigenvalues calculated with precision up to $10^{-15}$ we record it as the correct output.

\subsection{Sketch of algorithm}\label{sec:big_picture}
\noindent Frequently it is much easier to compute the largest eigenvalue of a matrix rather than the smallest eigenvalue, see for example \cite{KW}. At the same time, the smallest eigenvalue of an invertible matrix is the largest eigenvalue of its inverse. Therefore, it might prove advantages to treat the problem of calculating the smallest eigenvalue of a Hankel matrix as a problem of finding the inverse of the Hankel matrix and subsequently finding its largest eigenvalue.

The challenge with this approach is to compute the inverse of the Hankel matrix. Conveniently it is not necessary to compute the whole $H_N^{-1}$. Experimenting with the Hankel matrices in eq. (\ref{eq:matrix}) for relatively small size i.e. $N\sim100$, the authors have observed that largest eigenvalue of $H_N^{-1}$ is equal (to within $10^{-15}$ precision) to the largest eigenvalue of a $6\times6$ matrix constructed from the top-left entries of $H_N^{-1}$.

In other words, if one is interested in the largest eigenvalues of $H_N^{-1}$, one can discharge all entries of $H_N^{-1}$ except a small top-left $k \times k$ section, and compute the largest eigenvalues of that section. Moreover, $k$ does not need to be large. In the Table \ref{fig:table2} we found that $k=8$ was enough to achieve $10^{-12}$ precision for $N=4500$.

\subsection{{\it LDLT} algorithm}
\noindent In order to find the first $k \times k$ entries of $H_{N}^{-1}$ the authors have employed {\it LDLT} matrix decomposition.
In this subsection, we start by describing the sequential algorithm for computing the {\it LDLT} matrix decomposition.
In the next subsection, we follow by describing the parallelised version of the {\it LDLT} algorithm.

It is often considered the paradigm of numerical linear algebra that an algorithm must correspond to a matrix factorisation. The matrices we considered i.e. Hankel matrices are symmetric positive-semi-definite, therefore, it is natural to use Cholesky factorisation. In order to avoid the square root operations, we have to use the {\it LDLT} variant of Cholesky factorisation
\begin{align*}
 H = LDL^T
\end{align*}
where\footnote{
The matrix $H$ is positive semi-definite, as a result, the {\it LDLT} factorization will be numerically stable without pivoting.
Therefore, pivoting is not necessary.
Further, the authors suspect that for the case of Hankel matrices in eq. (\ref{eq:matrix}) any permutations would lead to slower execution and higher numerical errors.}
$L$ is a lower triangular matrix with ones along the diagonal and $D$ is a diagonal matrix.

\begin{table}
\begin{tabular}{p{13cm}}
\hline
\textbf{Algorithm 1: Serial LDLT code}\\
\hline
\begin{alltt}
  for (i=1 to N) \{
    {\it // this loop precomputes \(B\sb{i}\)}
    for (j=i+1 to N)
      \(B\sb{i}[j]=A\sb{i}[j]/A\sb{i}[i]\);

    {\it // these loops apply \(A\sb{i}\) and \(B\sb{i}\) to all entries to the right of \(A\sb{i}\)}
    for (j=i+1 to N)
      for (k=j to N)
        \(A\sb{j}[k]=A\sb{j}[k]-A\sb{i}[j]*B\sb{i}[k]\);
  \}
\end{alltt}
\\
\hline
\end{tabular}
\end{table}

In order to obtain the matrix factorisation in terms of $L$ and $D$, we employed the algorithm presented in the box named Algorithm 1.
Its outer loop is over the columns of $A$ from left to right.
The $i^{th}$ column of $A$ we will call $A_i$.
For each of those columns $A_i$, we first divide it by the diagonal entry i.e. $A_i[i]=A_{i\,i}$.
The result we call \(B\sb{i}\).
Subsequently we loop over all entries of the matrix to the right of $A_i$ and subtract
\begin{alltt}
 A[j][k] = A[j][k] - A[i][j]*B[i][k]
\end{alltt}
Notice that the term we subtract i.e. $A_i[j]*B_i[k]$ involves only the $A_i$ picked by the outer loop and $B_i$ we calculated at this step of the loop.

It is important to stress that the divisions involved in calculating the vector $B_{i}$
\begin{alltt}
 B[i][j] = A[i][j] / A[i][i]
\end{alltt}
are very expensive computationally. Therefore, significant amount of time is saved by precomputing it as indicated in the Algorithm 1 box rather than subtracting
\begin{alltt}
 A[j][k] = A[j][k] - A[i][j]*A[i][k]/A[i][i].
\end{alltt}

\begin{table}[!ht]
\centering
\caption{The balanced round robin assignment of columns to nodes.}
\label{table:assignment}
\large
\begin{tabular}{|lll|}
\hline
column $1$   & $\rightarrow$ & node $1$  \\
column $2$   & $\rightarrow$ & node $2$  \\
...          &               & ...       \\
column $n-1$ & $\rightarrow$ & node $n-1$\\
column $n$   & $\rightarrow$ & node $n$  \\
column $n+1$ & $\rightarrow$ & node $n$  \\
column $n+2$ & $\rightarrow$ & node $n-1$\\
...          &               & ...       \\
...          & $\rightarrow$ & node $2$  \\
...          & $\rightarrow$ & node $1$  \\
...          & $\rightarrow$ & node $2$  \\
...          &               & ...       \\
column $N-1$ & $\rightarrow$ & ...       \\
column $N$   & $\rightarrow$ & ...       \\
\hline
\end{tabular}
\end{table}

\subsubsection{Parallelising {\it LDLT}}
\noindent In order to parallelise the $LDLT$ decomposition on a computing cluster we followed the steps of \cite{ECW}. We assigned each column $A_i$ to one and only one nodes on the cluster. In order to assign similar amount of work to each node, we assigned columns to nodes using a (balanced) round robin approach. For $n$ nodes this approach assigns first $n$ columns to nodes $1$ up to $n$ and then next $n$ columns to nodes $n$ down to $1$, see table \ref{table:assignment}. This process gets repeated until each column is assigned to a node. This approach was found to produce satisfactory results for balancing the computational load and memory requirements on a homogeneous system in \cite{ECW} and no further improvements were introduced in that paper.

Splitting the columns between nodes enables the parallelisation of the ``for (j=i+1 to N)'' loop in the algorithm 1 over the nodes, see algorithm 2 box.

\begin{table}
\begin{tabular}{p{14.5cm}}
\hline
\textbf{Algorithm 2: The parallel version of the LDLT algorithm to be run by each node independently.}\\
\hline
\begin{alltt}
  Assign the columns to nodes in a balanced round robin fashion
  Broadcast the values needed to construct the initial matrix
  Compute \(B\sb{1}\)
  for (i=1 to N) \{
    if(column i+1 is assigned to this node) \{
      Apply the \(B\sb{i}\) to column \(A\sb{i+1}\)
      Compute \(B\sb{i+1}\)
      Initiate background transmit of \(A\sb{i+1}\) (and a part of \(B\sb{i+1}\))
      Apply \(A\sb{i}\) and \(B\sb{i}\) to all \(A\sb{i+2}\) ... \(A\sb{N}\) assigned to this node.
      Wait for transmit to complete
    \}
    else \{
      Initiate background receive of \(A\sb{i+1}\) (and a part of \(B\sb{i+1}\))
      Apply \(A\sb{i}\) and \(B\sb{i}\) to all \(A\sb{i+1}\) ... \(A\sb{N}\) assigned to this node.
      Wait for receive to complete
      Compute any missing \(B\sb{i+1}\) elements - discussed in detail below
    \}
  \}
\end{alltt}
\\
\hline
\end{tabular}
\end{table}

\newpage

\subsubsection{Communication between nodes}
\noindent In order to to distribute data between nodes and communicate between them we use OpenMPI library to implement MPI. As indicated in algorithm 2 box, initial distribution of the data between the nodes is done using MPI broadcasts. Subsequently at each loop iteration we broadcast \(A\sb{i+1}\) (and part of $B_{i+1}$) from the node that is assigned $A_{i+1}$ to all the other nodes.

We broadcast \(A\sb{i+1}\) (and part of $B_{i+1}$) as soon as they are calculated i.e. at $i^{th}$ step of the loop, see algorithm 2. Moreover, the communication is run in a separate thread that allows the transmission to overlap with the computation.

When sending $B_{i+1}$ we are faced with a choice:
\begin{itemize}
 \item we can broadcast the whole {\tt \(A\sb{i+1}\)} and make each node compute {\tt \(B\sb{i+1}\)} locally,
 \item or computing {\tt \(B\sb{i+1}\)} on the node that is assigned column {\tt \(A\sb{i+1}\)} and broadcast both {\tt \(A\sb{i+1}\)} and {\tt \(B\sb{i+1}\)} to the other nodes.
\end{itemize}
On one hand, the time needed to broadcast {\tt \(B\sb{i+1}\)} is significantly longer than the time needed to compute it. On the other hand, for small values of $i$ there is a significant amount of idle communication bandwidth. We use this bandwidth to broadcast both {\tt \(A\sb{i+1}\)} and {\tt \(B\sb{i+1}\)} in parallel to the computation. However, for smaller values of $i$ we limit the portion of {\tt \(B\sb{i+1}\)} we broadcast and let each node calculate the remaining part of {\tt \(B\sb{i+1}\)}. The final transmission algorithm is summarized in algorithm 3:

\begin{table}[H]
\begin{tabular}{p{14.5cm}}
\hline
\textbf{Algorithm 3:  The {\tt \(A\sb{i+1}\)} and {\tt \(B\sb{i+1}\)} transmission algorithm.}\\
\hline
\begin{alltt}
 serialize the \(A\sb{i+1}\) column
   send the \(A\sb{i+1}\) column
   break the \(B\sb{i+1}\) column into chunks of 100 values
   while (... there are more chunks ... and
          ... at least 8000 more multiplications to perform ...) \{
      send the next chunk of \(B\sb{i+1}\)
   \}
\end{alltt}
\\
\hline
\end{tabular}
\end{table}

The thresholds of 100 values and 8000 multiplications as presented in the algorithm 3 box were found empirically by the authors of \cite{ECW} ``to work well on a variety of problem sizes and cluster geometries''. This paper did not try to improve on those.

\subsubsection{Hybrid MPI and OpenMP parallelisation}
\noindent We have implemented the algorithm 2 using both MPI and OpenMP to manage the parallelism. In the preceding subsection we have just described how MPI is used to distribute columns between nodes. In contrast OpenMP is used to spread the computation assigned to a particular node across multiple cores within a node. In particular we have used OpenMP to parallelise the inner most loops associated with
\begin{itemize}
 \item compute $B_i$,
 \item apply the $A_i$ and $B_i$ to all \(A\sb{i+1}\) ... \(A\sb{N}\) assigned to this node,
\end{itemize}
see algorithm 2 box. The second of those loops is a loop over both rows and the columns assigned to a node. A primitive parallelisation would involve OpenMP parallelising one of those loops

\pagebreak

\begin{table}[H]
\begin{tabular}{p{14.5cm}}
\hline
\begin{alltt}
 for (j=i+1 to N) \{
   #pragma omp parallel for
   for (k=j to N)
     \(A\sb{j}[k]=A\sb{j}[k] - A\sb{i}[j]*B\sb{i}[k]\);
 \}
\end{alltt}
\\
\hline
\end{tabular}
\end{table}

However, it was noticed in \cite{ECW} that the implicit barrier at the end of the innermost loop resulted resulted in significant wast of computing time. Following \cite{ECW} we have implemented $A$ as a single index array which could be iterated over with a single loop

\begin{table}[H]
\begin{tabular}{p{14.5cm}}
\hline
\begin{alltt}
#pragma omp parallel for threadprivate(j,k) schedule(dynamic,chunk_size)
 for (index=lastIndex to firstIndex) \{
   j=... decode j from index ...
   k=... decode k from index ...
   \(V\sb{index}=V\sb{index} - A\sb{i}[j]*B\sb{i}[k]\);
 \}
\end{alltt}
\\
\hline
\end{tabular}
\end{table}

This loop is equivalent to the double loop above. Further improvements come from
\begin{itemize}
 \item dynamic scheduling - The computation time of each loop iteration varies considerably depending on the location in the matrix. Dynamic scheduling uses a queue to store block of work and assigns them the processor thread when they become available.
 \item $chunk\_size$ - Dynamic scheduling introduced a trade off. The larger the block size the more time wasted at the final barrier at the end of the loop, however, the smaller block size the larger overhead due to queue management overhead. We have used
 \begin{align}
   chunk\_size = \max \left(5, \dfrac{number\_of\_loop\_iterations}{200*number\_of\_OpenMP\_threads} \right).
 \end{align}
 \item We make the OpenMP loop run down through the indices. The underlying reason for going backwards is that the matrix entries become smaller toward the upper left corner of the matrix. This leads to the smaller execution blocks that lead to less time wasted at the final barrier at the end of the loop.
\end{itemize}

\subsection{Finding the inverse of $L$}\label{sec:invL}
\noindent In this subsection we will show that the {\it LDLT} matrix decomposition described in the previous section is the most computationally heavy step to find the first $k$ by $k$ entries of $H_{N}^{-1}$. Given $L$ and $D$ such that
\begin{align*}
 A = LDL^T
\end{align*}
in order to find the inverse
\begin{align*}
 H_{N}^{-1} = (L^{-1})^T  D^{-1}  L^{-1},
\end{align*}
or in index notation
\begin{align}\label{equ:index}
  \left(H_{N}^{-1}\right)_{jk} =  \sum_{l=1}^N L^{-1}_{lj} L^{-1}_{lk} / D_{ll}
\end{align}
we need to find $L^{-1}$. We find the inverse of $L$ by performing in-place Gauss-Jordan raw elimination. Let us first present the sequential version of the algorithm

\begin{table}[H]
\begin{tabular}{p{14.5cm}}
\hline
\textbf{Algorithm 3: The sequential version of in-place Gauss-Jordan raw elimination.}\\
\hline
\begin{alltt}
 for (i=0 to N-1)
 \{
    for (j=i+1 to N-1)
    \{
       f=A[j][i];
       A[j][i]=-f;

       for(k=0 to i-1)
          A[j][k]=A[j][k] - f*A[i][k];
    \}
 \}
\end{alltt}
\\
\hline
\end{tabular}
\end{table}

Moreover, we are only interested in the first $m$ by $m$ entries of $H_{N}^{-1}$, therefore, we are interested only the first $N$ by $m$ entries of $L^{-1}$. These can be found by a significantly faster algorithm, see algorithm 4 box below.

\begin{table}[H]
\begin{tabular}{p{14.5cm}}
\hline
\textbf{Algorithm 4: Sequential Gauss-Jordan raw elimination for the first $N$ by $m$ entries of $L^{-1}$.}\\
\hline
\begin{alltt}
 for (i=0 to N-1)
 \{
    for (j=i+1 to N-1)
    \{
       f=A[j][i];

       if(i < m)
          A[j][i]=-f;

       for(k=0 to min(m,i)-1)
          A[j][k]=A[j][k] - f*A[i][k];
    \}
 \}
\end{alltt}
\\
\hline
\end{tabular}
\end{table}

\subsubsection{Parallelisation of Gauss-Jordan elimination}
We parallelised the Gauss-Jordan elimination in a very similar way to the {\it LDLT} decomposition. We have distributing the rows of A over the nodes. We loop over rows of $A$ from top to bottom. At iteration number $j$ we broadcast the row number $i$ of $i$ to the rest of the nodes. After the transmission is finished, each node updates all rows assigned to it
\begin{alltt}
             A[k][j]=A[k][j] + A[i][j]*A[k][i];
\end{alltt}
We summarise the parallelised algorithm in Algorithm 5 box.

\begin{table}[H]
\begin{tabular}{p{14.5cm}}
\hline
\textbf{Algorithm 5:  The parallelised in-place Gauss-Jordan raw elimination for the first $N$ by $k$ entries of $L^{-1}$.}\\
\hline
\begin{alltt}
 Assign the columns to nodes in a balanced round robin fashion
 for (i=0 to N-1)
 \{
    n = the smallest of i and m
	
    if(row i is assigned to this node)
       Broadcast first n entries of i row of A to the other nodes
    else
       Receive first n entries of i row of A
    for (j=i+1 to N-1)
    \{
       if(row j is assigned to this node)
       \{
          f=A[j][i];
          if(i < k)
             A[j][i]=-f;

          for (k=0 to n-1)
             A[j][k]=A[j][k] - f*A[i][k];
       \}
    \}
 \}
\end{alltt}
\\
\hline
\end{tabular}
\end{table}

There are some important differences between how {\it LDLT} decomposition and Gauss-Jordan elimination parallelised
\begin{itemize}
 \item Unlike for {\it LDLT} decomposition, we use no OpenMP parallelisation. It would be natural to use OpenMP to parallelise the innermost for loop. However, we noticed that due to the upper limit ($=n$) the overhead was much larger than possible gains for small $m\sim10$).
 \item Due to small size of the broadcasted message which contains at maximum $m$ entries, we did not need a parallel communication and computation\footnote{One might suspect that an overhead from starting a parallel communication might actually result in a slower over all performance for small $m\sim10$.}. Table \ref{fig:table2} shows that during out tests the time spend communicating in the parallelised Gauss-Jordan was very short.
\end{itemize}

\subsection{Finding the inverse of Hankel matrix}\label{sec:invH}
\noindent Having finished the algorithm 4 we obtain the first $N$ by $k$ entries of $L^{-1}$ in place of the first $N$ by $k$ entries of $L$. Subsequently, we can perform the sum
\begin{align*}
  \left(H_{N}^{-1}\right)_{jk} =  \sum_{l=1}^N L^{-1}_{lj} L^{-1}_{lk} / D_{ll}
\end{align*}
to obtain the first $k\times k$ entries of $H_{N}^{-1}$. We perform the sum on the first node. The other nodes send the relevant entries of $L^{-1}$ as they become needed. The final values of $H_{N}^{-1}$ we cast into double precision floating point numbers. Finally, we use GNU Scientific Library to find the largest 2 eigenvalues of the $k\times k$ truncated $H_{N}^{-1}$.

When truncating $H_{N}^{-1}$ to the first $k\times k$ entries, we introduce a truncation error. In section \ref{sec:big_picture} we have argued that the truncation error is small, however, we would like to estimate it. For that purpose we calculate the largest 5 eigenvalues of $(k-1)\times (k-1)$ truncated $H_{N}^{-1}$. We then use the difference between the eigenvalues for $k$ truncation and $k-1$ truncation to estimate the truncation error, see table \ref{fig:table2}.

\section{Numerical results}\label{sec:numerics}
\noindent In section \ref{sec:methods}, we have presented a new numerical algorithm for computing the smallest eigenvalues of Hankel matrices, which we have implemented the algorithm in C programming language. In this section, we present the numerical results obtained using the new implementation.

\subsection{Computed eigenvalues}\label{sec:computed}
\noindent In this subsection, we present the computations we performed on High-Performance Computing Cluster of the University of Macau. Each computation was performed on 3 computing nodes, each node with 2 Intel® Xeon® E5-2690 v3 E5-2697 v2 CPUs and 64 GB of RAM. Each CPU had 12 cores/24 threads and 30 MB of cache. Therefore, in total, we had 36 cores and 192 GB of RAM available.

The numerical results obtained together with the timing profiles for the corresponding calculations we have presented in table \ref{fig:table2} on page \pageref{fig:table2}. The table columns correspond to
\begin{description}[noitemsep,topsep=0pt,parsep=3pt,partopsep=5pt]
 \item[N] is the rank of the matrix.
 \item[Required Precision] is the minimum precision of arithmetic operations $K$ to make the numerical error less than $10^{-15}$. We were increasing the value of $K$ in the steps of $1024$. The quoted minimum required precision $K$ was the smallest multiple of $1024$ such that the outputs of both $K$ and $K+1024$ matched to $10^{-15}$. For more details see section \ref{sec:precision} on page \pageref{sec:precision}.
 \item[No. of nodes] - the number of nodes used for the computation. For each computation, we have checked that during the computationally heavy part i.e. {\it LDLT} decomposition we have used 24 cores at (close to) 100\% at each node.
 \item[Total Wall time] is the total time taken by the computation as measured by ``a clock on the wall''.
 \item[{\it LDLT} time] is the time taken by the {\it LDLT} decomposition part of the algorithm.
 \item[Inversion of $L$ time] is divided into time spent performing arithmetic operations involved in finding the inverse of $L$, and time spent communicating between the nodes. We also quote the sum the above.
 \item[Transpose time] is the time taken to swap between entire columns being assigned to a particular node for the LDLT decomposition and entire rows being assigned to a particular node for the inversion of $L$.
 \item[Inverse of $H_N$] is the time taken by multiplication of two $L^{-1}$ and $D^{-1}$ matrices to find the truncated $H_N^{-1}$.
 \item[Truncation error] is the estimate of the error resulting from truncation of $H_N^{-1}$ to $k$ by $k$ size, see section \ref{sec:invH} for details.
\end{description}

All of the above times were measured on the host node i.e. the node that was recruiting other nodes to store assigned columns and to perform operations on the assigned columns see section \ref{sec:methods}.

There are a few important things to notice in table \ref{fig:table2}. First, the inversion of $L$ time constitutes only a few percents of the total wall times taken by the computation, which demonstrates that the {\it LDLT} decomposition is the computationally heavy part of the computation. Second, the time taken by matrix multiplication to find $H_N^{-1}$ constitutes an even smaller fraction of the total wall times taken by the computation, even though, it was performed on a single node. Third, the ``Send/Receive'' time is a small fraction of the ``Inversion of $L$ time'', therefore, parallelising the communication and arithmetics would not lead to significant speed-up.

\afterpage{%
    \clearpage
    \thispagestyle{empty}
     \begin{landscape}
     \captionof{table}{\footnotesize Numerically calculated smallest eigenvalues of Hankel matrix with $\beta=\frac{1}{2}$ together with the required precision and timing data of the corresponding computation. For an explanation of the column names see section \ref{sec:computed}.}
    \label{fig:table2}
     \setlength{\tabcolsep}{0.5pt}
     \renewcommand{\arraystretch}{1.5}
    {\scriptsize \begin{longtable}[l]{r|ccccc|ccc|c|cc|cc|cc}
    \hline\hline
      \thead{N}&\thead{\shortstack{Precision\\Required}}&\thead{\shortstack{No. of\\nodes}}&\thead{\shortstack{Total\\wall\\time}}&\thead{\shortstack{LDLT\\time}}&\thead{\shortstack{Transpose\\time}}&\multicolumn{3}{c}{\textbf{Inversion of $L$ time}} &\thead{\shortstack{Inverse\\of $H_{N}$}}&\multicolumn{6}{c}{\textbf{Calculated eigenvalues}}\\

    \thead{ }&\thead{ }&\thead{ }&\thead{ }&\thead{ }&\thead{ }&\multicolumn{3}{c}{ } &\thead{ }&\multicolumn{6}{c}{ }\\

        \thead{}&\thead{[bits]}&\thead{}&\thead{[s]}&\thead{[s]}&\thead{[s]}&\thead{Total}&\thead{Arithmetics}&\thead{\shortstack{Send/\\Receive}}&\thead{[s]}&\thead{Smallest}
        &\thead{\shortstack{Truncation\\error}} &\thead{$2^{nd}$ smallest}&\thead{\shortstack{Truncation\\error}}&\thead{$3^{nd}$ smallest}&\thead{\shortstack{Truncation\\error}}\\
        \hline
        \textbf{500}  & 1024 & 3 & 20.4  & 11.5  & 0.587 & 0.958 & 0.813 & 0.143 & 6.84 & 0.1204653471966412  & $-6.52\times 10^{-16}$ & 1.116696239796391  & $-2.75\times 10^{-14}$ & 33.53605844924584 & $-3.45\times 10^{-11}$ \\
                      &      &   &       &       &       &       &       &       &      &                     &                        &                    &                        &
                          &                         \\
        \textbf{1000} & 2048 & 3 & 136   & 109   & 3.83  & 11.7  & 11.5  & 1.23  & 9.17 & 0.08208748342129053 & $-5.57\times 10^{-15}$ & 0.8694471685364237 & $-2.51\times 10^{-13}$ & 16.74741576006559 & $-1.41\times 10^{-10}$ \\
        &      &   &       &       &       &       &       &       &      &                     &                        &                    &                        &
                          &                         \\
        \textbf{1500} & 3072 & 3 & 741   & 645   & 12.5  & 50.7  & 45.8  & 5.00  & 28.9 & 0.06529477501882298 & $-1.75\times 10^{-14}$ & 0.7587290286009394 & $-8.72\times 10^{-13}$ & 11.56571839375061 & $-2.98\times 10^{-10}$ \\
        &      &   &       &       &       &       &       &       &      &                     &                        &                    &                        &
                          &                         \\
        \textbf{2000} & 4096 & 3 & 2586  & 2374  & 28.3  & 142   & 129   & 12.9  & 33.4 & 0.05543072589537470 & $-3.68\times 10^{-14}$ & 0.6903595403385252 & $-2.02\times 10^{-12}$ & 9.032809963814945 & $-4.95\times 10^{-10}$ \\
        &      &   &       &       &       &       &       &       &      &                     &                        &                    &                        &
                          &                         \\
        \textbf{2500} & 4096 & 3 & 6717  & 6298  & 52.9  & 297   & 271   & 26.6  & 53.2 & 0.04878757749929328 & $-6.39\times 10^{-14}$ & 0.6418871023190091 & $-3.80\times 10^{-12}$ & 7.522048961034497 & $-7.23\times 10^{-10}$ \\
        &      &   &       &       &       &       &       &       &      &                     &                        &                    &                        &
                          &                         \\
        \textbf{3000} & 5120 & 3 & 15343 & 14578 & 94.9  & 581   & 530   & 51.0  & 61.6 & 0.04394036934849594 & $-9.83\times 10^{-14}$ & 0.6047656504596126 & $-6.28\times 10^{-12}$ & 6.513621908349736 & $-9.77\times 10^{-10}$ \\
        &      &   &       &       &       &       &       &       &      &                     &                        &                    &                        &
                          &                         \\
        \textbf{3500} & 6144 & 3 & 29001 & 27818 & 140   & 940   & 933   & 6.99  & 71.4 & 0.04021149503682476 & $-1.40\times 10^{-13}$ & 0.5749057442617896 & $-9.52\times 10^{-12}$ & 5.789838905207254 & $-1.25\times 10^{-09}$ \\
        &      &   &       &       &       &       &       &       &      &                     &                        &                    &                        &
                          &                         \\
        \textbf{4000} & 7168 & 3 & 56445 & 54426 & 227   & 1645  & 1502  & 142   & 86.3 & 0.03723304780176154 & $-1.88\times 10^{-13}$ & 0.5500577166937035 & $-1.36\times 10^{-11}$ & 5.243332063875005 & $-1.55\times 10^{-09}$ \\
        &      &   &       &       &       &       &       &       &      &                     &                        &                    &                        &
                          &                         \\
        \textbf{4500} & 7168 & 3 & 94242 & 91303 & 325   & 2436  & 2217  & 218   & 94.1 & 0.03478615399760864 & $-2.43\times 10^{-13}$ & 0.5288610646385768 & $-1.84\times 10^{-11}$ & 4.814940754432488 & $-1.87\times 10^{-09}$ \\
        \hline\hline
        \end{longtable}}
    \end{landscape}
    \clearpage
}

\subsection{Timing: the new algorithm against Secant algorithm}\label{sec:compare}
\noindent The authors of \cite{ECW} compared the Secant algorithm with a number of classical eigenvalue algorithms including Householder, Jacobi, Lanczos on the task of finding the smallest eigenvalue of Hankel matrices. It was found that for the case of large and extremely ill-conditioned matrices Secant algorithm proved to be much more efficient. In this paper, we try to improve on the Secant algorithm with the new algorithm presented in section \ref{sec:methods}. Here we would like to compare the timing of the new algorithm implementation against the implementation of Secant algorithm provided by \cite{ECW}.

In table \ref{table:compare} we have compared the time needed to compute the smallest eigenvalue of $N$ by $N$ Hankel matrix from eq. (\ref{eq:matrix}) to $10^{-15}$ numerical accuracy with the two algorithms for different values of $N$. For each $N$ we have scanned over precision $K$, defined in eq. (\ref{eq:precision}), and chose the minimum value that produced the desired accuracy and called it the ``required precision''. The values of precision $K$ we scanned over were 256, 388, 512, 768, 1024, 1536, 2048, 3072, 4096 bits. The computation was performed on Thinkpad T480 with Intel® Core™ i5-8250U Processor using close to 100\% of all 8 threads for both implementations.

The table contains the number of iterations the Secant algorithm needs to converge on the zero of the characteristic polynomial of the matrix for each $N$, let us call it $n_{iterations}$. The total wall time taken by the Secant implementation was composed of $n_{iterations}$ {\it LDLT} decompositions plus the initial decomposition. Therefore
\begin{align*}
 t_{\text{Wall Secant}} = (n_{iterations}+1) \times t_{\text{Average iteration Secant}}
\end{align*}
Looking at table \ref{table:compare}, we can notice two improvements over the Secant algorithm:
\begin{itemize}
 \item The computation time required by the new implementation is very close to the time required for the Secant implementation to complete an average iteration. Therefore, the new implementation is $\sim(n_{iterations}+1)$ times faster. This is to be expected as the new algorithm performs the computationally heavy {\it LDLT} decomposition only once, whilst the Secant algorithm does that $(n_{iterations}+1)$ number of times.
 \item We notice that the precision $K$ required to compute the smallest eigenvalue with the desired precision is typically significantly smaller for the new implementation. This can be attributed to the fact that the Secant algorithm finds a zero of the characteristic polynomial whose values are typically astronomically large, and therefore, finding the zero to satisfactory precision is harder than finding a ``satisfactory'' {\it LDLT} decomposition. This leads to small but significant speed gains for larger matrices.
\end{itemize}

\begin{table}[H]
\centering
\caption{Comparison between the new algorithm and Secant algorithm used in \cite{ECW}. For an explanation of the column names see section \ref{sec:computed}.}
\label{table:compare}
\setlength{\tabcolsep}{2.3pt}
\renewcommand{\arraystretch}{1.6}
{\normalsize\begin{longtable}{c|cccc||cccc}
\hline\hline
\thead{} & \multicolumn{4}{c}{\bf Secant} & \multicolumn{4}{c}{\bf New algorithm} \\

\thead{N}& \thead{\shortstack{Precision\\Required}}&\thead{\shortstack{Wall\\time}}&\thead{\shortstack{No. of\\iterations}}&\thead{\shortstack{Average\\iteration}}
& \thead{\shortstack{Precision\\Required}}&\thead{\shortstack{Wall\\time}}&\thead{\shortstack{LDLT\\time}}&\thead{\shortstack{Inv. L\\time}}\\
\thead{}&\thead{[bits]}&\thead{[s]}&\thead{}&\thead{[s]}& \thead{[bits]}&\thead{[s]}&\thead{[s]}&\thead{[s]}\\
\hline
200 & 512 & 230.6 & 8 & 25.6 & 388 & 22.4 & 22.2 & 0.037 \\
400 & 1024 & 505.4 & 8 & 56.15 & 758 & 55.1 & 54.2 & 0.527 \\
600 & 1024 & 1652.9 & 11 & 137.7 & 1024 & 134.9 & 131.8 & 2.24 \\
800 & 2048 & 2395.7 & 8 & 266.2 & 1536 & 376.4 & 367.9 & 6.79 \\
1000 & 3072 & 8542.4 & 7 & 1067.8 & 2048 & 954.6 & 935.5 & 15.7 \\
\hline\hline
\end{longtable}}
\end{table}

\subsection{Scaling with number of nodes}\label{sec:nodes}
\noindent Authors of \cite{ECW} extensively discuss how and present evidence that their implementation of Secant algorithm scales well with the number of nodes. In this section, we would like to check that this statement extends to the new algorithm.

We have timed our implementation of the new algorithm for a different number of nodes. We have used Google Cloud Platform to create a cluster of 6 ``f1-micro'' virtual machines\footnote{These machines have $2.5$ GHz CPU and $0.6$ GB of memory.} and run our implementation on a various number of nodes to check the scaling of the computation time with the number of nodes. The results are presented in table \ref{table:nodes}. The row names follow the same convention as table \ref{fig:table2} with one difference; we also quote the total CPU time which is the sum of time consumed by all of the CPUs utilized by the program. For this setup we estimate the CPU time using
\begin{align*}
 t_{CPU} \simeq n_{nodes} \times t_{Wall}.
\end{align*}

We can see that the total CPU time only slowly increases with the number of nodes, which leads to the total wall time steadily decreasing with the number of nodes. We can see that we achieve parallelisation over the nodes of the cluster, and it scales well with the number of recruited nodes.

\begin{table}[H]
\centering
\caption{Timing test of the scalability of the algorithm with the number of nodes participating in the calculation. For an explanation of the column names see section \ref{sec:nodes}.}
\label{table:nodes}
\begin{tabular}{lllllll}
\textbf{Number of nodes}    & \textbf{}                                  & \textbf{2} & \textbf{3} & \textbf{4} & \textbf{5} & \textbf{6} \\ \hline
\textbf{Total wall time}    & \multicolumn{1}{l|}{\textbf{}}             & 29.3       & 20.1       & 16.5       & 13.8       & 12.6       \\
\textbf{Total CPU time}     & \multicolumn{1}{l|}{\textbf{}}             & 58.6       & 60.4       & 66.1       & 69.2       & 75.5       \\
\textbf{LDLT time}          & \multicolumn{1}{l|}{\textbf{Total wall}}   & 27.4       & 18.4       & 14.7       & 12.3       & 11         \\
\textbf{}                   & \multicolumn{1}{l|}{\textbf{Total CPU}}    & 54.9       & 55.1       & 58.9       & 61.7       & 66.3       \\
\textbf{Inverse of $L$ time} & \multicolumn{1}{l|}{\textbf{Total wall}}  & 0.672      & 0.465      & 0.366      & 0.301      & 0.263      \\
\textbf{}                   & \multicolumn{1}{l|}{\textbf{Total CPU}}    & 1.34       & 1.40       & 1.46       & 1.51       & 1.58       \\
\textbf{}                   & \multicolumn{1}{l|}{\textbf{Arithmetics}}  & 0.608      & 0.424      & 0.311      & 0.257      & 0.205      \\
\textbf{}                   & \multicolumn{1}{l|}{\textbf{Send/Receive}} & 0.0641     & 0.041      & 0.055      & 0.044      & 0.058      \\
\textbf{Inverse $H_N$ time} & \multicolumn{1}{l|}{\textbf{}}             & 0.22       & 0.158      & 0.175      & 0.169      & 0.176      \\
\textbf{Eigenvalue time}    & \multicolumn{1}{l|}{\textbf{}}             & 0.004      & 0.034      & 0.034      & 0.002      & 0.002
\end{tabular}
\end{table}

\section{Leading exponent determination}\label{sec:exponent}
In this section, we use the numerically calculated eigenvalues for $\beta=\frac{1}{2}$ from section \ref{sec:computed} to compare it to eq. (\ref{eq:leading_order}) and extract the leading exponent introduced in section \ref{sec:sad}.

The authors of \cite{ECW} have also compared eigenvalues they have numerically calculated for $\beta=\frac{1}{2}$ to eq. (\ref{eq:leading_order_old}). Their approach was to directly compare the numerically calculated eigenvalue and the corresponding value predicted by eq. (\ref{eq:leading_order_old}). They have observed a significant difference\footnote{The difference is significantly smaller when we compare the numerics to eq. (\ref{eq:leading_order}) instead of eq. (\ref{eq:leading_order_old}).} even for relatively large $N$ e.g. $>25\%$ for $N=1500$. These differences are due to the subleading in $N$ order term not captured by eq. (\ref{eq:leading_order}).

In contrast, we would like to concentrate on the scaling of the smallest eigenvalue with $N$
\begin{align*}
 \lambda_1(N) \propto N^{-\frac{2}{\pi}} \sqrt{\log{N}}
\end{align*}
and emphasise the algebraic factor $N^{\frac{-2}{\pi}}$ in our analysis. For that purpose we recast eq. (\ref{eq:leading_order}) which reads
\begin{align}\label{eq:leading_order2}
 \lambda_1(N)               &\simeq \dfrac{8\pi\sqrt{\log{N}}}
                                           {(4\pi N {\rm e})^{\frac{2}{\pi}}}.
\end{align}
as
\begin{align}\label{eq:linear}
 \log\left[\dfrac{8\pi}{\lambda_1(N)}\sqrt{\log N}\right] \simeq \dfrac{2}{\pi} \log(4\pi  N{\rm e}).
\end{align}
We notice that if we plot $\log\left[\frac{8\pi}{\lambda_1(N)}\sqrt{\log N}\right]$ against $\log(4\pi N{\rm e})$ for large $N$ we should find a straight line with the gradient equal to the leading exponent i.e. $\frac{2}{\pi}$. For smaller $N$ we would expect some deviation from a straight line due to the subleading terms not captured by eq. (\ref{eq:leading_order2}).

This picture emphasises the leading exponent. Moreover, as we have noted in section \ref{sec:sad}, the subleading terms do not affect the overall $N^{-\frac{2}{\pi}}$ factor, and they are suppressed by a double log-function. Therefore, we would expect the leading exponent determination to be less affected by the subleading terms than the direct comparison performed by \cite{ECW}.

\subsection{Linear fit}
\begin{figure}[H]
 \centering
 \includegraphics[width=0.95\textwidth]{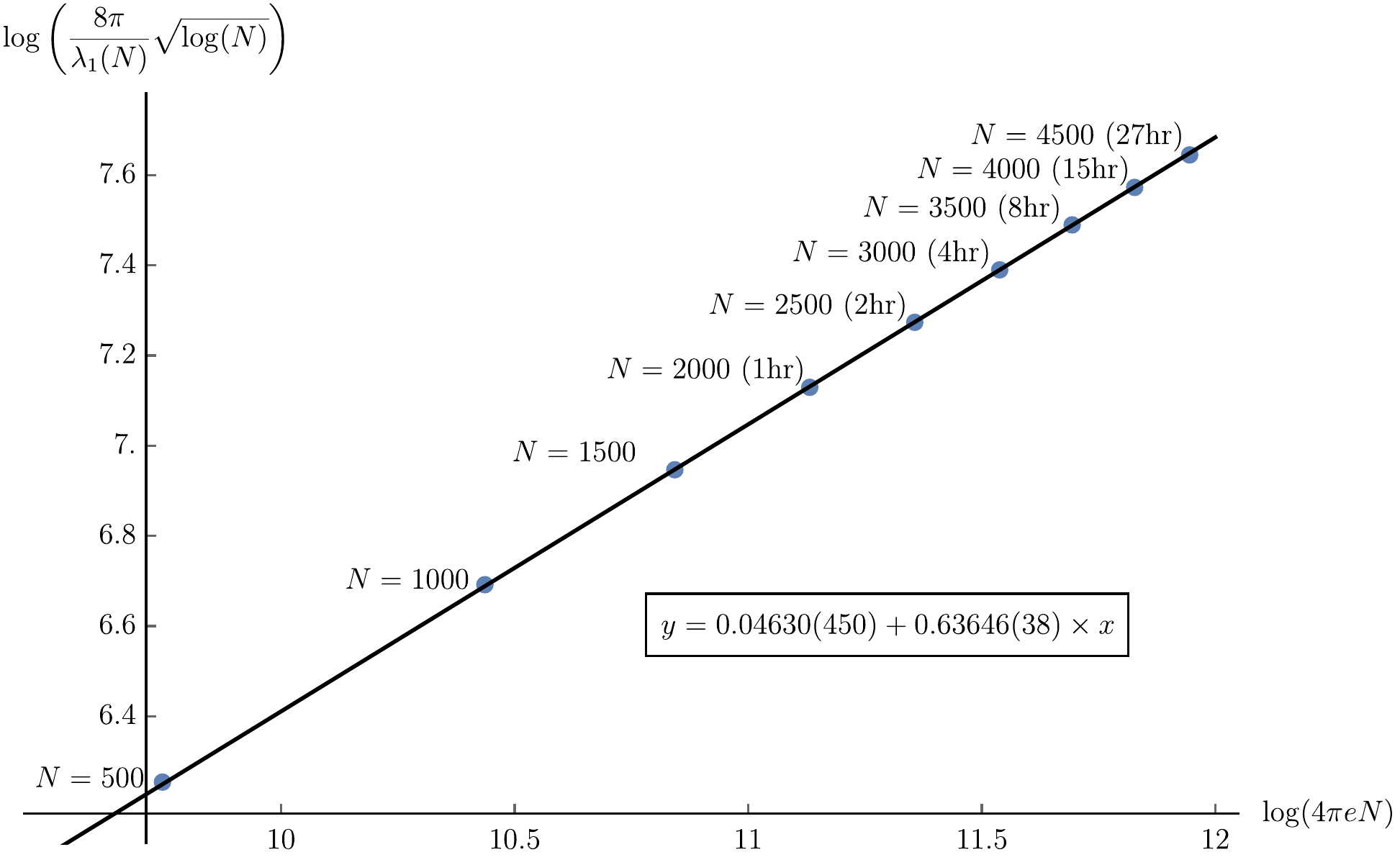}
 \caption{Plot of the numerically calculated eigenvalues from table \ref{fig:table2} with a fitted linear model.}
 \label{fig:plot}
\end{figure}
\begin{figure}[H]
 \centering
 \includegraphics[width=0.85\textwidth]{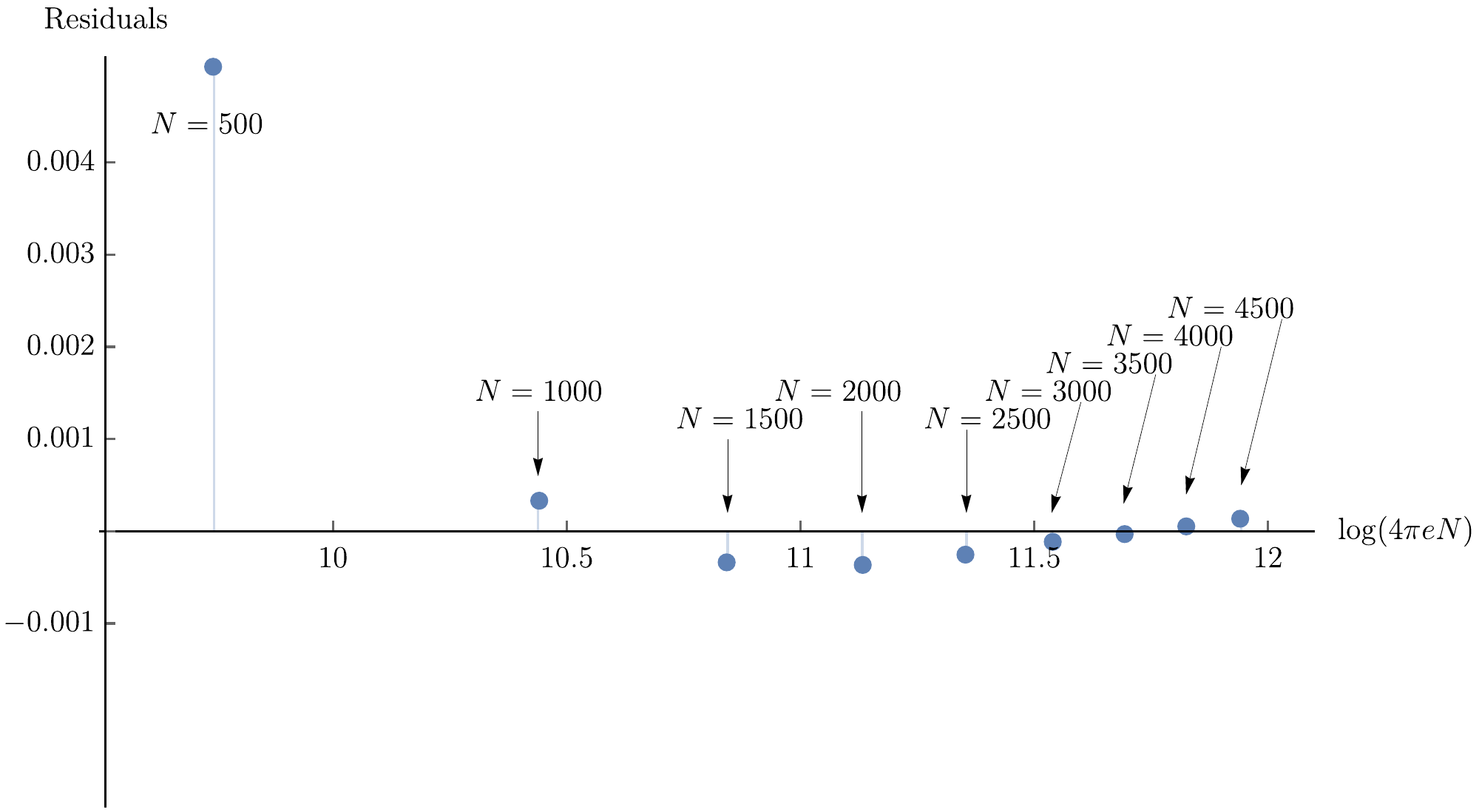}
 \caption{The residuals of the linear fit in fig. \ref{fig:plot}.}
 \label{fig:residues}
\end{figure}

We plotted the eigenvalues presented in table \ref{fig:table2} on a x-y plot where
\begin{align}
 x &\equiv \log(4\pi N{\rm e}),\\
 y &\equiv \log\left(\dfrac{8\pi}{\lambda_1(N)}\sqrt{\log N}\right),
\end{align}
in figure \ref{fig:plot}. Subsequently, we have fit the data points with a linear model. We have included a constant term in the linear model to allow for the differences for smaller values of $N$. We weighted each point proportionally to $N^2$ to put higher weight on the larger values of $N$ and found the best fit is
\begin{align*}
 y = 0.04630 + 0.63646 \times x
\end{align*}
with adjusted $\text{R-squared}=0.999997$, and residuals presented in figure \ref{fig:residues}.
We can see that the residuals tend to decrease as $N$ increases. The gradient value is only $0.025\%$ away from the predicted value $\frac{2}{\pi}\simeq0.63662$, i.e.
\begin{equation*}
\frac{0.63662-0.0.63646}{0.63662}\simeq0.00025
\end{equation*}
Moreover,  the $95\%$ confidence interval\footnote{We should mention that the points in figure \ref{fig:plot} have negligible (numerical and truncation) error bars associated with them, but since we are fitting the model parameters to the data points those model parameters acquire statistical errors.} for the gradient $[0.63510, 0.63781]$ contains the $\frac{2}{\pi}$. Therefore, we conclude that the smallest eigenvalue of the Hankel matrices in eq. (\ref{eq:matrix}) at the critical point $\beta=\frac{1}{2}$ decays algebraically $\lambda_1(N) \sim N^{-\frac{2}{\pi}}$ for large $N$, and the numerical results are a strong evidence that eq. (\ref{eq:poly}) and eq. (\ref{eq:leading_order}) are true asymptotically for large $N$.

\section{Discussion}
\noindent In this paper, we have introduced a novel algorithm for finding the smallest eigenvalues of a matrix that is particularly suited for Hankel matrices generated by the probability density function $w(x) = {\rm e}^{-x^\beta}$. We have developed an implementation of the algorithm that runs on a high performance computing cluster making use of both shared memory parallelisation within a node and a distributed memory parallelisation between nodes of the cluster. We have shown that the new algorithm improves on Secant algorithm previously considered to be the state of the art algorithm to study of Hankel matrices \cite{ECW}. Moreover, we show that it scales well with the number of participating nodes.

We have employed the novel algorithm to find the smallest eigenvalues of Hankel matrices at the critical point $\beta=\frac{1}{2}$ for matrices up to $N=4500$. We have plotted the numerical results on a custom log-log plot and determined that at the leading exponent is $-\frac{2}{\pi}$, therefore
\begin{align*}
 \lambda_1(N) \simeq N^{-\frac{2}{\pi}},
\end{align*}
where $\lambda_1(N)$ is the smallest eigenvalue of the Hankel matrix in eq. (\ref{eq:matrix}) at $\beta=\frac{1}{2}$ critical point.

\subsection{Significance of the $\beta=\frac{1}{2}$ result}
\noindent Berg, Chen and Ismail \cite{BCI} proved that the minimum eigenvalues of Hankel matrices generated by probability density function $w(x)$
\begin{equation*}
\left(H_{N}\right)_{ij}=\int_{0}^{\infty}x^{i+j}w(x)dx,~~0\leq i,j\leq N-1,
\end{equation*}
is bounded away from zero for all $N$ including $N$ equals infinity if and only if the associated moment problem has more than one solutions. For the weight $w(x)={\rm e}^{-x^{\beta}},~x\geq0,~\beta>0$ it can be shown that the smallest eigenvalue tends to zero exponentially fast in $N$ if $\beta>\frac{1}{2}$, a critical point, where the moment problem is at the verge of being indeterminate, the smallest eigenvalue tends to zero algebraically in $N$. For $0<\beta<\frac{1}{2}$, the smallest eigenvalue tends to a strictly positive constant whose value depending on $\beta$ is not yet found.

\subsection{Preconditioner}
\noindent The authors suspect that there might be space for further improvements in the design of the numerical algorithm to calculate the smallest eigenvalue. The Hankel matrix we consider has particularly natural and effective preconditioner and post-conditioner that make its condition number much more tame. One can consider
\begin{align*}
  \left(\tilde{H}_{N}\right)_{ij} &= \dfrac{\Gamma\left(\dfrac{i+j-1}{\beta}\right)}{\sqrt{\Gamma\left(\dfrac{2i-1}{\beta}\right)\Gamma\left(\dfrac{2j-1}{\beta}\right)}}, & &\text{for } i,j=1,2,3...
\end{align*}
which has a nice property that $(\tilde{H}_{N})_{ij} = 1$ for $i=j$ and a significantly smaller condition number. Therefore, it might be significantly easier to find the determinant of $(\tilde{H}_{N})_{ij}$ than the original $(H_{N})_{ij}$. However, $(\tilde{H}_{N})_{ij}$ is just the original Hankel matrix pre-multiplied and post-multiplied by two diagonal matrices. Therefore, the determinant of the original $(H_{N})_{ij}$ is just the determinant of $(\tilde{H}_{N})_{ij}$ multiplied by the diagonal terms of pre and post-multiplier. This sounds like a powerful idea to improve Secant algorithm, which relays on a repeated evaluation of (modified) determinants of Hankel matrix. However, we were surprised to find out that the pre-conditioner and post-conditioner did not bring any efficiency improvements to the Secant algorithm. This might be partially explained by the fact that $(\tilde{H}_{N})_{ij}$ still has a very large condition number though much smaller than $(H_{N})_{ij}$. We feel that this situation deserves further clarification, and we leave it as a future direction.

\section{Acknowledgements}
\noindent{\color{black}Y. Chen, J. Sikorowski and M. Zhu would like to thank the Science and Technology Development Fund of the Macau SAR for generous support in providing FDCT 130/2014/A3 and FDCT 023/2017/A1. We would also like to thank the University of Macau for generous support via MYRG 2014-00011 FST, MYRG 2014-00004 FST} and {\color{black}MYRG 2018-00125 FST.}


\begin{thebibliography}{99}
\bibitem{C16}
Akhiezer NI. \emph{The Classical Moment Problem and Some Related Questions in Analysis}. Edinburgh: Oliver and Boyd; 1965.
\bibitem{C24}
Beckermann B. The condition number of real Vandermonde, Krylov and positive definite Hankel matrices. \emph{Numer Math.} 2000; 85: 553-557.
\bibitem{BCI}
Berg C, Chen Y, Ismail MEH. Small eigenvalues of large Hankel matrices: The indeterminate case. Math. Scand. 2002; 91: 67-81.
\bibitem{C14}
Berg C, Szwarc R. The smallest eigenvalue of Hankel matrices. \emph{Constr Approx}. 2011; 34: 107-133.
\bibitem{C23}
Chen Y, Ismail MEH. Thermodynamic relations the Hermitian matrix ensembles. \emph{J Phys A: Math Gen}. 1997; 30: 6633-6654.
\bibitem{CL1}
Chen Y, Lawrence N. Small eigenvalues of large Hankel matrices. J. Phys. A: Math. Gen. 1999; 32: 7305-7315.
\bibitem{C10}
Chen Y, Lubinsky DS. Smallest eigenvalues of Hankel matrices for exponential weights. \emph{J Math Anal Appl}. 2004; 293: 476-495.
\bibitem{C27}
De Bruijn NG. \emph{Asymptotic Methods in Analysis}. New York: Interscience; 1958.
\bibitem{ECW}
Emmart N, Chen, Y, Weems C. Computing the smallest eigenvalue of large ill-conditioned Hankel matrices. Commun. Comput. Phys. 2015; 18: 104-124.
\bibitem{zzz}
Ismail MEH. \emph{Classical and Quantum Orthogonal Polynomials in One Variable}, Encyclopedia of Mathematics
and its Applications, Vol. 98. Cambridge, UK: Cambridge University Press; 2005.
\bibitem{KW}
Kuczy\'{n}ski J, Wo\'{z}niakowski H. Estimating the largest eigenvalue by the power and lanczos algorithms with a random start. SIAM Journal on Matrix Analysis and Applications. 1992; 13: 1094-1122.
\bibitem{new11}
Krein MG, Nudelman AA. \emph{Markov Moment Problems and Extremal Problems}. Providence, RI: American Mathematical Society; 1977.
\bibitem{C25}
Lubinsky DS. Condition numbers of Hankel matrices for exponential weights. \emph{J Math Anal Appl}. 2006; 314: 266-285.
\bibitem{C17}
Mehta ML. \emph{Random Matrices, Third Edition}. Singapore: Elsevier (Singapore) Pte Ltd.; 2006.
\bibitem{C3}
Szeg\"{o} G. On some Hermitian forms associated with two given curves of the complex plane. \emph{Trans Amer Math Soc}. 1936; 40: 450-461.
\bibitem{C5}
Todd J. Contributions to the solution of systems of linear equations and the determination of eigenvalues. \emph{Nat Bur Standards Appl Math Ser}. 1959; 39: 109-116.
\bibitem{C6}
Widom H, Wilf HS. Small eigenvalues of large Hankel matrices. \emph{Proc Amer Math Soc}. 1966; 17: 338-344.
\bibitem{new33}
Widom H, Wilf HS. Errata: Small Eigenvalues of Large Hankel Matrices. \emph{Proc Amer Math Soc}. 1968; 19: 1508.
\bibitem{C22}
Zhu M, Chen Y, Emmart N, Weems C. The smallest eigenvalue of large Hankel matrices. \emph{Applied Mathematics and Computation}. 2018; 334: 375-387.
\end{thebibliography}
\end{document}